\documentclass{article}

  \usepackage{amsmath}
\usepackage{layout}
\usepackage{graphicx}
 \newcommand\toe{\stackrel{e}{\to}}
\newcommand\tov{\stackrel{v}{\to}}

\begin{document}

\title{BOUNDS  ON SOME EDGE FOLKMAN NUMBERS \footnote{Partially
supported by the scientific research fund of St. Kl. Ohridski  Sofia
University under contract No 75/2009 } }

\maketitle

\author{Nikolay R. Kolev}

\begin{abstract}
{  
   The edge Folkman numbers are defined by the equality
  $$ F_e(a_1, \ldots , a_r ;q) = \min\{ | V(G)| : G \toe  (a_1, \ldots ,
   a_r; q) \text{ and } cl(G)<q  \}.
  $$
  Lin proved that if $R(a_1, a_2)=R(a_1-1, a_2)+R(a_1, a_2-1)$ then  $F_e(a_1 ,
   a_2 ; R(a_1, a_2))=R(a_1, a_2)+2$, where $R(a_1, a_2)$ is the Ramsey number
   showing that $ K_{R(a_1, a_2)-3} +C_5 \toe (a_1, a_2)$
  where $C_5$ is the simple cycle on 5 vertices. We prove some upper
   bounds on edge Folkman numbers for which
   $R(a_1, a_2)< R(a_1-1, a_2)+R(a_1, a_2-1)$
  and we cite some lower bounds.
 }
\end{abstract}

\section{1. Introduction}
\small{
Only finite non-oriented graphs
without multiple edges and loops are considered. We call a
$p$-clique of the graph $G$ a set of $p$ vertices each two of which
are adjacent. The largest positive integer $p$ such that $G$
contains a $p$-clique is denoted by cl(G). A set of vertices of the
graph $G$ none two of which are adjacent is called an independent
set. The largest positive integer $p$ such that $G$ contains an
independent set on $p$ vertices is called the independence number of
the graph $G$ and is denoted by $ \alpha(G).$ In this paper we shall
also use the following notations:
\begin{itemize}
\item  $V(G)$ is the vertex set of the graph $G$;
 \item $E(G)$ is the edge set of the graph $G$;
\item   $N(v)$, $v\in V(G)$ is the set of all vertices of
$G$ adjacent to $v$;
 \item  $G[V]$, $V \subseteq V(G)$ is the subgraph of $G$ induced by $V$;
  \item    $K_n$ is the complete graph on $n$ vertices;
  \item $\overline{G}$ is  the complementary graph of $G $.
  \end{itemize}
   }

\small{ Let $G_1$ and $G_2$ be two graphs without common vertices.
We denote by $G_1+G_2$ the graph $G$ for which $V(G)=V(G_1)\cup
V(G_2)$ and $E(G)=E(G_1)\cup E(G_2)\cup E'$ where $E'=\{ xy:x\in
V(G_1), y\in V(G_2)\}$. It is clear that

\begin{equation}
 cl(G_1 + G_2) = cl(G_1) + cl(G_2).
\end{equation}
}


\small{ \textbf{ Definition 1.} Let $a_1, \ldots , a_r$ be positive
integers. The symbol $G \tov (a_1, \ldots , a_r)$ means that in
every $r$-coloring of $V(G)$ there is a monochromatic $a_i$-clique
in the $i$-th color for some $i \in \{1, \ldots, r\}$.
\medskip


\textbf{ Definition 2.} Let $a_1, \ldots , a_r$ be positive
integers. We say that an $r$-coloring of $E(G)$ is $(a_1, \ldots ,
a_r)$-free if for each $i=1, \dots , r$ there is no monochromatic
$a_i$-clique in the $i$-th color. The symbol $G \toe (a_1, \ldots ,
a_r)$ means that there is no $(a_1, \ldots , a_r)$-free coloring of
$E(G).$

 The smallest positive integer $n$ for which $K_n \toe (a_1, \ldots ,
   a_r)$ is  called a Ramsey number and is denoted by $R (a_1, \ldots ,
   a_r).$  Note that the Ramsey number $R(a_1, a_2)$ can be
   interpreted as the smallest positive integer $n$ such that for
   every $n$-vertex graph $G$ either $cl(G)\geq a_1 $ or $\alpha(G)\geq a_2.
   $ The existence of such numbers was proved by Ramsey in [18]. We
   shall use only the values
   $R(3,3)=6$ and $R(3,4)=9,$ [3].

   The edge Folkman numbers are defined by the equality
 $$ F_e(a_1, \ldots , a_r ;q) = \min\{ | V(G)| : G \toe  (a_1, \ldots ,
    a_r) \text{ and } cl(G)<q  \}.
 $$
It is clear that $G\toe (a_1,\ldots ,a_r)$ implies $cl(G)\geq
\max\{a_1,\ldots , a_r\}$. There exists a graph $G$ such that $G
\toe (a_1, \ldots , a_r) $ and $cl(G) = \max\{a_1,\ldots , a_r\}$.
In the case $r=2$ this was proved in [1] and in the general case in
[16]. Therefore
  $$\label{1} F_e(a_1 ,\ldots , a_r ;q) \mbox{ exists
    if and only if }
  q> \mbox{max} \{ a_1,\ldots , a_r \}.$$

It follows from the definition of $R(a_1,\dots,a_r)$ that
\[
F_e(a_1,\dots,a_r;q)=R(a_1,\dots,a_r) \text{ if }q>R(a_1,\dots,a_r).
\]

The smaller the value of $q$ in comparison to $R(a_1,\dots,a_r)$ the
more difficult the problem of computing the number
$F_e(a_1,\dots,a_r;q).$ Only ten edge Folkman numbers that are not
Ramsey numbers are known. For the results see the papers: [2],[4],
[6], [7], [9], [10], [11], [13].}


\section{2. Upper bounds on the edge Folkman numbers}

\small{ We
obtained the following result:

\textbf{Theorem 1} \textit{Let}  $a $ \textit{and} $\alpha$ \textit{
be nonnegative   integers}. \textit{  Let us
      denote } $R=R(3,a).$ \textit{ Let }
      $R= R(3,a -1) + a- \alpha $  \textit{ and }
     $ R  -3a + \alpha +5 \geq R(3,a-2)$,  $ a \geq 4$ .
    \textit{If there
    exists a graph} $U$ \textit{with the properties}
    $$ cl(U) =a -1  $$
    $$ U \tov (a-1, a-2) $$
  $$ U \tov (  \underbrace {a-3, \ldots ,a-3}_{a-2 \mbox{  times}},3) .$$
 \textit{ Then}
$$ F_ e (3,a; R-a+ \alpha + 4) \leq    R -2a+ \alpha +4   + |V(U)|.  $$

\textit{Proof}  Consider the graph $G = K_{R -2a+ \alpha +4} +
  U  $ where $U$ is the graph from the statement of the theorem. We
will  prove that $G \toe (3,a)$ and thus the theorem will be proved
(here we use (1) to compute that $cl(G)= R -2a+ \alpha +4 +cl(U)=
R -a +\alpha +3 $).

Assume the opposite: that there exists a $(3,a)$-free 2-coloring of
$E(G) $
\begin{equation}
 E(G)=E_1\cup E_2, \qquad E_1\cap E_2=\emptyset.
\end{equation}
We shall call the edges in $E_1$   blue and the edges in $E_2$  red.

We define for an arbitrary vertex $v\in V(G)$ and index $ i=1,2 :$
\begin{align*}
N_i(v)&=\{x\in N(v)\mid [v,x]\in E_i\},
 \\
G_i(v)&=G[N_i(v)]\\
A_i(v) &= N_i(v) \cap V(U)
\end{align*}

Let $H$ be a subgraph of $G$. We say that $H$ is a monochromatic
subgraph in the blue-red coloring (2) if $E(H)\subseteq E_1$ or
$E(H)\subseteq E_2$. If $E(H)\subseteq E_1$ we say that $H$ is a
blue subgraph and if $E(H)\subseteq E_2$ we say that $H$ is a red
subgraph.

It follows from the assumption that the coloring (2) is
$(3,a)$-free that
 \begin{equation}
      cl(G_1(v)) \leq a-1
 \text{ and }
   cl(G_2(v)) \leq  R(3, a-1) -1
   \text{ for each }  v \in V(G)
 \end{equation}
 Indeed, assume that $cl(G_1(v)) \geq a.$ Then there must be no blue
 edge connecting any two of the vertices in $  G_1(v)  $ because
 otherwise this blue edge together with the vertex $v$ would give a
 blue triangle. As we assumed $cl(G_1(v)) \geq a$ then we have a red
$a$-clique. Analogously assume $cl(G_2(v)) \geq  R(3,a-1)$.
 Then we have either a blue 3-clique or a red $(a-1)$-clique in $G_2(v).$ If
 we have a blue 3-clique in $G_2(v)$ then we are through. If we have
 a red $(a-1)$-clique   then this $(a-1)$-clique together with the vertex $v$ gives a
 red $a$-clique. Thus (3) is proved.

 We shall prove that
\begin{equation}
 cl(G[A_1(v)])+cl(G[A_2(v)])\leq 2a-5
 \text{ for each }
  v \in V(K_{R -2a+ \alpha +4})
 \end{equation}
Assume that (4) is not true, i.e. that there exists a vertex $v
\in V(K_{R -2a+ \alpha +4})$ such that

$$cl(G[A_1(v)])+cl(G[A_2(v)])\geq 2a-4.$$

Then as there  are $R -2a+ \alpha +3$ more vertices in $  V(K_{R
-2a+ \alpha +4})$ with the exception of $v,$ it follows that
$$cl(G_1(v)) +cl(G_2(v)) \geq R -2a+ \alpha +3+2a-4 = R + \alpha -1  = R(3,a-1) +a -1              $$
(here we use the statement of the theorem that $R(3,a)=R(3,a-1)+a -
\alpha$.)

This contradicts
(3). Thus (4) is proved.

Now we shall prove that
    \begin{equation}
 cl(G[A_1(v)])=a-1 \text{ or } cl(G[A_2(v)])=a-1
 \text{ for   each }  v \in    V(K_{R-2a+ \alpha +4})
 \end{equation}

Assume that (5) is not true. Then we obtain from $cl(U) =
a-1$  that

\begin{equation}
  cl(G[A_1(v)])\leq a-2 \text{ and } cl(G[A_2(v)]) \leq a-2~\text{for some}~v \in   V(K_{R-2a+ \alpha +4}).
 \end{equation}

It follows from the statement of the theorem that $ U \tov (a-1,
a-2) $ that in every 2-coloring of $V(U),$ in which there are no
$(a-1)$-cliques in none of the two colors then there are $(a-2)$-cliques
in the both colors.Hence $A_1(v)$
and $A_2(v)$ contain $(a-2)$-clique.  Therefore the inequalities in (6) are in fact
equalities, which contradicts (4). Thus (5) is proved. Let us
note that it follows from (5) and (4) that for each $ v \in
V(K_{R-2a+ \alpha +4})$

\begin{equation}
 \text { if } cl(G[A_2(v)])=a-1 \text { then }  cl(G[A_1(v)]) \leq a-4.
 \end{equation}
 If we assume that there are $a$ vertices $v$ in
$V(K_{R-2a+ \alpha +4})$ with the property  $cl(G[A_1(v)]) = a-1 $
then it follows from (3) that there are only red edges between
these $a$ vertices and hence this is a red $a$-clique.  Therefore
having in mind the statement of the theorem    there are at least
$R-2a+ \alpha +4 -(a-1) \geq R(3,a-2) $ vertices $v$ in $V(K_{R-2a+
\alpha +4})$ with the property $cl(G[A_1(v)])\leq a-2 .$ Let us
denote the set of these vertices by $S$. It follows by (5) that
$cl(G[A_2(v)]) = a-1 $ for each $v \in S.$ Then we obtain from (7)
that
\begin{equation}
  cl(G[A_1(v)])\leq a-4 \mbox{ for each } v \in S .
\end{equation}
As there are no blue 3-cliques among the vertices in $S$
  there is a red $(a-2)$-clique among them.
Let us denote the vertices of this $(a-2)$-clique by $w_i,$ $i \in
1,\ldots ,a-2.$ Let us partition the vertices of $U$ in $a-1$ colors
in the following way:

$$V_1 =  A_1(w_1)$$

and for $j \in 2, \ldots, a-2 $
$$V_j = A_1(w_j)  \backslash V_1  \cup   \ldots \cup V_{j-1}$$

$$V_{a-1}=   V(U) \backslash V_1  \cup   \ldots \cup V_{a-2}. $$

 According to (8) it is impossible    $ V_j$, $j \in 1,\ldots ,a-2$ to contain  a
$(a-3)$-clique because according to its definition  $V_j$ is a subset
of $A_1(w_j)$. But we know from the statement of the theorem
 that $ U \tov (  \underbrace {a-3, \ldots ,a-3}_{a-2 \mbox{  times}},3) .$
  Hence $V_{a-1} $ contains a 3-clique. Then this 3-clique must contain a
red edge (otherwise it is a blue 3-clique and we are through). But
it follows from the definition of  $V_{a-1} $ that all the edges
between the vertices in $V_{a-1} $ and the vertices  $w_1, \ldots ,
w_{a-2} $ are red. Then the red edge in $V_{a-1} $ and the red
$(a-2)$-clique  $w_1, \ldots , w_{a-2} $ form a red $a-$clique. Thus
the theorem is proved.

\medskip

If we put   $\alpha=0,$  $a=5,$ $U=Q,$
 where $Q$ denotes the graph
whose complementary is the graph
     given here  we obtain the following Kolev`s result
     $$F_e(3,5;13) \leq 21, [5].$$

  \begin{figure}
\centering
 \includegraphics{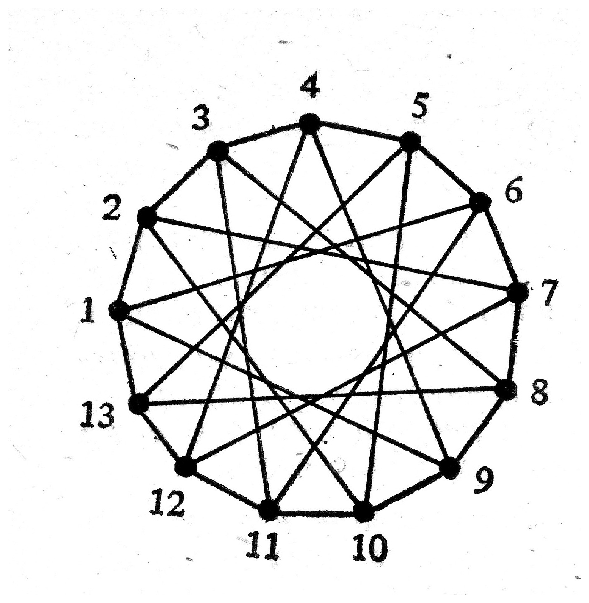}
\caption{Graph $\overline{Q}$}
\end{figure}

What is novel in this result is that for the first time the graph
$Q$ is used in the theory of edge Folkman numbers.

             We also prove the following result which unfortunately is not a particular case of
             \textbf{Theorem1}:

\textbf{Theorem 2}     $F_e(4,4;17) \leq 25.$

\textit{Proof }  Consider the graph $ G=K_{12}+Q $ and we shall show
that $G \toe (4,4).$

Assume the opposite: that there exists a $(4,4)$-free 2-coloring of
$E(G) $
\begin{equation}
 E(G)=E_1\cup E_2, \qquad E_1\cap E_2=\emptyset.
\end{equation}
We shall call the edges in $E_1$   blue and the edges in $E_2$  red.

We define for an arbitrary vertex $v\in V(G)$ and index $ i=1,2 :$
\begin{align*}
N_i(v)&=\{x\in N(v)\mid [v,x]\in E_i\},
 \\
G_i(v)&=G[N_i(v)]\\
A_i(v) &= N_i(v) \cap V(Q)
\end{align*}

  We say that $H$ is a monochromatic
subgraph in the blue-red coloring (9) if $E(H)\subseteq E_1$ or
$E(H)\subseteq E_2$. If $E(H)\subseteq E_1$ we say that $H$ is a
blue subgraph and if $E(H)\subseteq E_2$ we say that $H$ is a red
subgraph.

It follows from the assumption that the coloring (9) is (4,4)-free
and $R(3,4)=9$ that
 \begin{equation}
      cl(G_1(v)) \leq 8
 \text{ and }
   cl(G_2(v)) \leq 8
   \text{ for each }  v \in V(G)
 \end{equation}

   Assume $cl(G_2(v)) \geq 9$. Since $R(3,4) =9$,
 then we have either a blue 3-clique or a red 4-clique in $G_1(v).$ If
 we have a blue 3-clique in $G_1(v)$ then this blue 3-clique together with the vertex $v$
forms a blue 4-clique.
  If we have
 a red 4-clique   then we are through.   Analogously we disprove  $cl(G_1(v)) \geq 9$. Thus (10) is proved.

 We shall prove that
\begin{equation}
 cl(G[A_1(v)])+cl(G[A_2(v)])\leq 5
 \text{ for  each }
  v \in V(K_{12})
 \end{equation}
Assume that (11) is not true, i.e. that there exists a vertex $v
\in V(K_{12})$ such that

$$cl(G[A_1(v)])+cl(G[A_2(v)])\geq 6.$$

Then as there  are eleven more vertices in $  V(K_{12})$ with the
exception of $v,$ it follows that
$$cl(G_1(v)) +cl(G_2(v)) \geq 17.             $$
It follows from the pigeonhole principle that either $cl(G_1(v))\geq
9$ or $cl(G_2(v))\geq 9,$ which contradicts (10). Thus (11) is
proved.

It follows from $Q \toe (3,4)$ and $cl(Q) =4$ that
    \begin{equation}
 cl(G[A_1(v)])=4 \text{ or } cl(G[A_2(v)])=4
 \text{ for each }  v \in V(K_{12}).
 \end{equation}
It follows from (11) and (12) that

\begin{equation}
 \text { if }cl(G[A_1(v)])=4 \text { then }  cl(G[A_2(v)]) \leq 1
  \text{ for each }  v \in V(K_{12}).
  \end{equation}

 Assume that there are two vertices
$ a,b \in V(K_{12})$ such that the edge $ab$ is blue and $
cl(G[A_1(a)]) = cl(G[A_1(b)]) =4.$ Then it follows from (13) that
$ cl(G[A_2(a)])\leq 1$ and $ cl(G[A_2(b)]) \leq 1.$ Let us consider
the following 3-coloring of $V(Q)$:
$$ V_1=A_2(a)$$
$$V_2=A_2(b) \backslash A_1(a)$$
$$ V_3 = V(Q) \backslash (V_1 \cup V_2).$$
It is clear from the definition of $V_1$, $V_2$, $V_3$ that $V_3$
consists of such vertices that are connected with both $a$ and $b$
with blue edges. It is clear from  $cl(G[A_2(a)])\leq 1$ and $
cl(G[A_2(b)]) \leq 1 $ that $V_1$ and $V_2$ contain no 2-cliques.
But we know  from [8] that $Q \tov (2,2,4)$ and hence $V_3$ contains
a 4-clique. If this 4-clique has a blue edge then then this  blue
edge together with the edge $ab$ is a blue 4-clique. If this
4-clique does not have a blue edge then it is a red 4-clique- a
contradiction. Thus we obtained that

\begin{equation}
\text{if }cl(G[A_1(a)]) = cl(G[A_1(b)]) =4 \text{ then the edge } ab
\text{ is not blue, } a,b \in K_{12}.
 \end{equation}
Analogously we obtain

\begin{equation}
\text{if }cl(G[A_2(a)]) = cl(G[A_2(b)]) =4 \text{ then the edge } ab
\text{ is not red, } a,b \in K_{12}.
 \end{equation}

We have from (12) that $ cl(G[A_1(v)])=4 \text{ or }
cl(G[A_2(v)])=4$ for each $v \in K_{12}.$ For a fixed  $v \in
K_{12}$ we may assume without loss of generality that $
cl(G[A_1(v)])=4$. Assume that for this $v$ there exist three
distinct vertices $b_1, b_2, b_3 \in K_{12}$ such that the edges $v
b_i$, $i \in \{1, 2,3 \}$ are blue. Now it follows from
(13), (14),(15) and the fact that  the edges $v b_i$, $i \in
\{1, 2,3 \}$ are blue that

\begin{equation}
cl(G[A_2(b_i)])=4 \text { and }  cl(G[A_1(b_i)]) \leq 1
 \end{equation}
 for  $i \in \{1, 2,3 \}.$
Now it follows from (15) that the edges $b_i b_j$ are blue
for $i, j \in \{1, 2,3 \}.$ Hence $v b_1 b_2 b_3$ is a blue
$4-$clique and the theorem is proved.

Thus we obtained that  there are no more than two blue edges from
the vertex $v$. Therefore there are at  least 9 vertices connected
to $v$ with red edges which contradicts (10). We proved that $G \toe (4,4)$.
As $cl(G)=4$ then we obtain that $F_e(4,4;17) \leq 25$, which we wanted to prove.

\textbf{
     Remark } So far it was known  by [1]
      that $F_e(4,4;17) <  \infty$.

The last known  lower bound of the edge Foklman number $F_e(4,4;17)$ is
$F_e(4,4;17) \geq 22$ - see [15].


\section{References}

\leftskip 2pc
\parindent -2pc

\dotfill    

\small {[1]Folkman, J. Graphs with monochromatic complete subgraphs
in every edge coloring. \emph{SIAM J. Appl.~Math.} \textbf{18},
1970, 19--24.  [2]  Graham R.L. On edgewise 2-colored graphs with
monochromatic triangles containing no complete hexagon \emph{J.
Comb. theory} \textbf{4}, 1968, 300. [3]Greenwood, R., A. Gleason.
Combinatorial relation and chromatic graphs. \emph{Canad.~J. Math.,
\textbf{7}, 1955, 1--7.} [4]  Lin S. On Ramsey numbers and
$K_r-$coloring of graphs \emph{J. Comb. theory Ser B} \textbf{12},
1972, 82-92. [5]  Kolev N. New upper  bound for the edge Folkman
number $F_e(3,5;13).$ \emph{Serdica math. J.,   \textbf{34}
2008,783-790.} [6]Kolev, N., N. Nenov. An example of 16-vertex
Folkman edge (3,4)-graph without 8-cliques. \emph{ Annuaire Univ.
Sofia Fac. Math. Inform., \textbf{98}, 2008  101- 126, see
$http://arxiv.org/PS_cache/math/pdf/0602/0602249v1.pdf.$} [7] Kolev
N., N. Nenov. The Folkman number $F_e(3,4;8)$ is equal to 16. \emph{
Compt. rend.~Acad.~bulg.~Sci., \textbf{59},No 1, 2006,   25--30.}
[8]Nenov N.On the vertex Folkman number F(3,4). \emph{Compt.
rend.~Acad.~bulg.~Sci., \textbf{54}, 2, 2001, 23-26.} [9] Nenov N.
On an assumption of Lin about Ramsey-Graham-Spencer numbers.
\emph{(Russian)Compt. rend.~Acad.~bulg.~Sci., \textbf{33}, 9, 1980,
1171-1174.} [10] Nenov N. Generalization of a  certain theorem of
Greenwood and Gleason on three-color colorings of the edges of a
complete graph with 17 vertices. \emph{(Russian)Compt.
rend.~Acad.~bulg.~Sci., \textbf{34},  1981, 1209-1212.}
 [11] Nenov N. Lower estimates for
some constants related to Ramsey graphs.\emph{(Russian) Annuaire
Univ. Sofia Fac. Math. Inform., \textbf{75},  1984, 27-38.}
 [12]Nenov N. On the (3,4)-Ramsey graphs without 9-cliques.
\emph{(Russian) Univ. Sofia Fac. Math. Inform., \textbf{85},
 1991 ,71-81 .}
 [13] Nenov N. An example of 15-vertex Ramsey (3,3)-graph
with clique number 4. \emph{(Russian)Compt. rend. ~Acad. ~bulg.
~Sci., \textbf{34}, 1981, 1487-1489.}
 [14] Nenov N.On the Zykov numbers and some of its applications to  Ramsey theory.
 \emph{(Russian)Compt. rend. ~Acad. ~bulgSerdica Bulg. math. Publ.{9},161-167, (1983).}
  [15] Nenov N.  On the vertex Folkman numbers$F_v(2, \ldots, 2;q).$
  \emph{Serdica Math. J. 35, 2009  251-272.}
[16]Nesetril J,  Rodl V. The Ramsey property for graphs with
forbidden complete subgraphs.\emph{J. Comb Th. \textbf{B20}, 1976,
243-249.}
 [17] Piwakowski K., Radziszowski S., Urbanski S.
 Computation of the Folkman number $ F_e (3,3;5).$
 \emph{J. Graph  theory  } \textbf{32}, 1999, 41-49.
[18]Ramsey P. On a problem of formal  logic\emph{   Proc. London
Math. Soc., 30, 1930,  264-286.} }
 \vskip 0.5cm


{\it Faculty of mathematics and informatics
\hfill  Received xx.xx.199x  

Sofia University

Sofia 1164, BULGARIA

E-MAIL: nickyxy@fmi.uni-sofia.bg }

\end{document}